\documentclass[twoside,12pt]{amsart}
\pdfoutput=1 
\usepackage{comment}
\usepackage{mathrsfs,amsmath,amssymb,amsthm,epsfig}
\usepackage[all]{xy}
\usepackage{graphicx}
\def\rev{\usepackage[active]{srcltx}}\rev
\usepackage[colorlinks,plainpages,backref]{hyperref}
\usepackage{pdf14} 

\newcommand{\escala}[1]{\setlength{\unitlength}{#1cm}}
\newcommand{\red}[1]{{\color{red}{#1}}} 
\newcommand{\blu}[1]{{\color{blue}{#1}}}
\newcommand{ \af}[1]{\ensuremath{{\mathbb A}^{#1}}} 
\newcommand{ \pd}[1]{\ensuremath{\check{\mathbb P}^{#1}}} 
\newcommand{ \cl}[1]{\ensuremath{\mathcal{#1}}} 
\newcommand{\bsm}{\left(\begin{smallmatrix}}
\newcommand{\esm}{\end{smallmatrix}\right)}
\newcommand{\ba}[1]{\begin{array}{#1}}\newcommand{\ea}{\end{array}}
\newcommand{\na}[1]{\noalign{\vskip#1pt}}
\newcommand{\bb}[1]{\ensuremath{\mathbb{#1}}} 
\newcommand{\mb}[1]{\ensuremath{\mathbf{#1}}}
\newcommand{\wed}[1]{\ensuremath{\stackrel{#1}\wedge}}
\newcommand{\bmp}{\begin{minipage}}
\newcommand{\emp}{\end{minipage}}
\newcommand{\bc}{\begin{center}}
\newcommand{\ec}{\end{center}}
\newcommand{\be}{\begin{equation}}
\newcommand{\ee}{\end{equation}}
\newcommand{\nn}{\ensuremath{\ne}}
\newcommand{\bi}{\begin{itemize}}
\newcommand{\ei}{\end{itemize}}
\newcommand{\ov}[1]{\ensuremath{\overline{#1}}}
\newcommand{\ds}[1]{\displaystyle{#1}}
\newcommand{\us}{\ensuremath{^{\star}}}
\newcommand{\ls}{\ensuremath{_{\star}}}
\newcommand{\ve}{\ensuremath{{^\vee }}}
\newcommand{\lar}{\ensuremath{\longrightarrow}}

\newcommand{\id}[1]{\ensuremath{\langle{#1}\rangle}}
\newcommand{\wt}[1]{\ensuremath{\widetilde{#1}}}
\newcommand{\wh}[1]{\ensuremath{{\widehat{#1}}}}

\newcommand{\s}[1]{\ensuremath{{\operatorname S}_{#1}}}
\newcommand{\G}{{\bb G}}
\newcommand{\pol}{polynomial}
\newcommand{\lra}{\ensuremath{\leftrightarrow}}
\newcommand{\surj}{\ensuremath{
\xymatrix
{\,\ar@{->>}[r]&\,}}}

\newcommand{\F}[1]{\ensuremath{\bb F\!\operatorname{ol}_{#1}}}
\newcommand{\C}{\ensuremath{\bb C}}
\newcommand{\D}{\ensuremath{\bb D}}
\newcommand{\st}{\ensuremath{\,|\,}}
\newcommand{\uu}[1]{\ensuremath{\underline{#1}}}
\newcommand{\p}[1]{\ensuremath{\mathbb P^{#1}}} 
\newcommand{\vez}{\ensuremath{{\times}}}
\newcommand{\inv}{\ensuremath{{}^{\!-1}}}
\newcommand{\ra}{\ensuremath{\,\rightarrow\,}}
\newcommand{\ie}{{\em i.e., }}
\newcommand{\sF}{\ensuremath{\mathscr{F}}}

\begin{document}
\title[{Enumerative geometry of 
Legendrian foliations}]
{Enumerative geometry of 
\\Legendrian foliations: a Tale of Contact} 
\author[Mauricio Corr\^ea]{Mauricio Corr\^ea}
\address{\sc Mauricio Corr\^ea\\
\newline
\indent
Universit\`a degli Studi di Bari, 
Via E. Orabona 4, I-70125, Bari, Italy
\newline
\indent
UFMG\\
Avenida Ant\^onio Carlos, 6627\\
30161-970 Belo Horizonte\\ Brazil
}
\email{mauriciojr@ufmg.br, mauricio.barros@uniba.it}

\author[Israel Vainsencher]{Israel Vainsencher$^\star$}
\thanks{$^\star$Partially supported by CNPq
303755/2019-9 and FAPEMIG  RED-00133-21.}
\address{\sc Israel Vainsencher\\
\newline
\indent
UFMG\\
Avenida Ant\^onio Carlos, 6627\\
30161-970 Belo Horizonte\\ Brazil}
\email{israel@mat.ufmg.br}


\begin{abstract}
 \noindent
A contact distribution on \p3 is defined by the
1-form $\omega:=x_2dx_1-x_1dx_2+x_4dx_3-x_3dx_4$,
up to a change of projective coordinates.  The family
of contact distributions
 is parameterized by the complement of the
Pfaff-Pl\"ucker quadric in the projective 5-space
of antisymmetric 4\vez4 matrices. A foliation of
dimension 1 and degree $d$ is specified by a
polynomial vector field $\varphi:=\sum
p_i\partial_{x_i}, \,p_i$ homogeneous of degree
$d$. 
The foliation is called Legendrian if tangent to some
distribution of contact. Our goal is to give
formulas for the dimensions and degrees of the
varieties of Legendrian foliations, and of the
varieties of 
foliations tangent to a pencil of planes.
\begin{flushright}\medskip
To Alberto Collino, 
{\em in memoriam}.
\end{flushright}
\end{abstract}

\maketitle

\section{Introduction}

A holomorphic \emph{contact structure} on a
complex manifold $X$ is a codimension one
maximally non-integrable holomorphic sub-bundle
$\mathscr{F} \subset TX$.  The pair
$(X,\sF)$ is called a contact complex manifold.  Such
structures have appeared in a link
between Riemannian and Algebraic Geometries. 
For instance, Salamon \cite{Salamon} and LeBrun
\cite{Lebrun} proved that K\"ahler-Einstein
Fano contact manifolds are twistor spaces of
positive quaternion-K\"ahler manifolds. See also
Beauville\,\cite{Beauville}. 
Demailly\,\cite{Demailly} and
Kebekus\,{\em et\,al.}\,\cite{KPSW} have
 shown that projective contact manifolds
 are either Fano with $b_2=1$ or
a projectivized tangent bundle of a projective
manifold.  We refer the reader to
Druel\,\cite{Druel,Druel2} and Ye\,\cite{Ye} for
further classification results in low dimensions
and toric manifolds.  
In
\cite{Kebekus} Kebekus shows that if $(X,\sF)$ is
a Fano contact manifold with $b_2=1$ which is
not a projective space, then $X$ is covered by
lines tangent to the contact distribution
$\sF$. Such lines are examples of 
\emph{Legendrian curves}, $f:C\to
(X,\sF)$ with $f_*TC\subset \sF|_C$.
Enumerative aspects of Legendrian curves of 
contact  in 
$\mathbb{P}^{2n+1}$ and their
moduli spaces have been studied  in
\cite{Amorim,Bryant,Kalinin,LM,LV,Muratore}.

We investigate \emph{Legendrian foliations} on
\p3. These are foliations of dimension
one\,(cf.\,\ref{spfolns}) whose leaves are tangent
to a contact structure (\ref{foltgd}). They appear
in the classification of foliations of low degree
for which the singular schemes are of pure
dimension one, cf.\,\cite{CJM}. 

We present a  description of
the spaces of Legendrian foliations enabling us to
determine a formula for their degrees,
in the spirit of \cite{FV, LeV, RV,VK}. 

A {\em distribution of degree m and codimension one}
on \p3 is defined by a 1-form $\sum {a_i}
dx_i$ where the $a_i$ are homogeneous polynomials
of degree $m+1$ such that $\sum {a_i}x_i=0.$ We
assume henceforth $m=0$.  As syzygies-trained minds
will recognize, this entails the vector of
coefficients $(a_1,\dots,a_4)$ is a linear
combination $\sum a_{ij}\kappa_{ij}$ of  the six basic
Koszul  relations $\kappa_{12}:=
(x_2,-x_1,0,0),\dots,\kappa_{34}:=(0,0,x_4,-x_3)$. 
Each distribution of degree 0 is specified by a
point in the 
projectivization, \p5,  of the space of anti-symmetric
4\vez4 matrices, cf.\,\S\ref{syz}.
Matrices of maximal rank correspond to the {\it
 distributions of contact}. These distributions 
form the
 Zariski open subset of indecomposable 1-forms, 
complement of the Pfaff-Pl\"ucker quadric 
$\bb G\subset \p5$.
 The anti-symmetric matrices of rank 2 \ie
 elements in \bb G,
correspond to pencils of planes; this has lead us to
study also foliations tangent to a pencil of
planes,\,cf.\,\ref{foltgpencil}.

Grosso modo, a {\em foliation of dimension one} is
a (polynomial) recipe to draw a direction at each
point. Precisely, for each integer $d\geq0$ we
call 
\\\centerline{$\F d:=\p{}(H^0(\p3,T\p3(d-1)))$,}
 the space
of foliations of dimension 1 and degree
$d$. Each element $\varphi \in\F d$ is defined by
a polynomial vector field $\varphi:=\sum
p_i\partial_{x_i},$ with $p_i$ homogeneous of degree
$d$: think of the line joining
$P,\varphi(P)$ cf.\,\ref{spfolns}.  
It corresponds to a map of vector bundles,
$\cl O_{\p3}(1-d)\ra T\p3$. This defines  a
subspace of dimension one of the tangent space
$T_P\p3$ whenever $\varphi(P)\nn0$, \ie
$P$ is not a singular point of the foliation $\varphi$.

A foliation is called {\it Legendrian} whenever
tangent to some (variable) distribution of
contact, say $\omega$. This means the line
specified by the foliation $\varphi$ at a general
point is contained in the plane assigned by
$\omega,$ that is, $\omega\cdot\varphi=0$. For
fixed $\omega$, the condition is linear on
$\varphi$, thus defining a subspace
$\omega_d^\perp:=\{\varphi\in \F d\,
|\,\omega\cdot\varphi=0\}$. The main technical
difficulty stems from the fact that the dimension
of $\omega_d^\perp$ jumps as the 1-form of contact
$\omega\in\p5\setminus\bb G$ specializes to a
decomposable 1-form $\omega_0\in\bb G$ (\ie a
pencil of planes,\,see\,\S\ref{nml}). The
analogous question for tangency to other families
of distributions seems out of reach for our
toolbox.

Our main results  are  formulas for 
the degrees of the subvarieties of Legendrian
foliations (and friends: the subvarieties of
foliations tangent to a pencil of planes), 
see\,(\ref{athus1}),\,(\ref{athusbis}).
It turns out that the 
answers are given by polynomial functions on the
degree of the foliation.  This fits a trend of the
last few decades in Enumerative Geometry. We think
of it as a sort of Schubert Calculus programme for
exploring the geometry of parameter spaces of
foliations by curves. 

It is reminiscent of Alberto Collino's interest on
Enumerative Geometry of the Hilbert schemes of
lines (resp. conics) contained in certain
hypersurfaces, cf.\, \cite{col1,col2}.
The 2nd author misses the enlightening company of
Alberto pacing  the corridors at MIT back in the
early 70's.

\section{A couple of pictures for starters}
\label{fig1}
Each line in \p3 gives a recipe to draw a plane
through a general point: just take their linear span.
$$
\escala{1.5}\begin{picture}(.1,3.9)(6.,1.29)

\ifpdf
\put(2.976,1.6) {
\includegraphics[scale=0.67]
{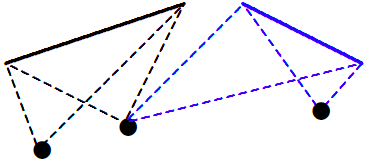}}
\fi
 \put(2.4,4.5) {A distribution of planes}

\put(6.492,4.5){\blu{Another distribution of planes}}

\put(4.08,2.15){\red{$\linethickness{2pt}\line(1,0){2} $}}

\end{picture}
$$

Their \red{intersections} yield a foliation of
dimension 1.  The distributions depicted above are
pencils of planes; they are defined by integrable
1-forms such as $\omega:=udv-vdu,$ with $u,v$
homogeneous polynomials of degree 1. The axis
$u=v=0$ is the singular locus of this
distribution.

\subsection{Another recipe for choosing 
a plane 
through each point in 3-space} 

Given $ (x_0,y_0,z_0)\in\af3$,
{\em le-voil\`a}: 
\bc $y-y_{0 }= z_0x   -x_0z .  
$\ec

We have a bijection, \bc
{\em point} \,$(x_0,y_0,z_0)\leftrightarrow
\text{\em plane }\,
z_0x -y   -x_0z  +y_{0 }$.
\ec

This bijection
extends to the  so-called
{\em null correlation} on \p3,
$$
\p3\,\ni\,\left[x_0,y_0,z_0,w_0\right]\,\leftrightarrow
\, z_0x -w_0y   -x_0z  +y_{0 }w\,\in\, \pd3,
$$
which is related  to a nice, classical 
construction sketched below for the reader's benefit.

\subsection{Revisit twisted   cubics}
\label{fig2}
Eleven out of ten geometers cherish the twisted
cubic as the most beloved rational curve in space.
It is at center stage in the above null
correlation and contact distribution.
 
For each point \ ${P}$ \ there are three
osculating planes drawn from \ {$P$} \ to the
cubic curve.  Take the connecting plane \,{$p$}\,
of the {osculating points}: 
${P} \leftrightarrow\red{ p=\id{O_1,O_2,O_3 } }.$

$
\escala{1.5}\begin{picture}(2,5.2)(1.5, 1.21)
\ifpdf
\put(1.3,1.7) {
\includegraphics[scale=0.7]
{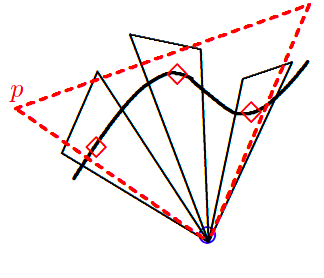}}
\fi
\put( 5.671,2.725){\fbox{\!
\bmp{ 7.57cm}
Exercise:   the plane   {$ p$}   contains \\
  the point 
{$P$},\,cf.\,\cite{Steinmetz},\cite[p.\,312]{Semple-Kneebone}.
\emp\!\hskip-2cm}}
\put(4.9 ,1.85){$P$}
\put(3.2 ,3.75){$\red{O_1}$}
\put(4.2 ,4.605){$\red{O_2}$}
\put(5.45 ,4.7){$\red{O_3}$}
\end{picture}$

Conversely, to every plane corresponds the point
$P$ of intersection of the osculating planes of
the twisted cubic curve, drawn in the 3 points of
intersection of the {plane} with the  curve.

\section{Distributions, 1-forms}

Let \s k\ denote the space of homogeneous \pol s
of degree $k$.
A holomorphic distribution of degree $m$
of planes in \p3 is defined by a projective 1-form
\be\label{dist}
\omega:=a_1dx_1 + \cdots + a_4dx_4,\,
a_i\in{\s{m+1}},
\ee
such that 
$
\sum a_ix_i =  0 . $ 
For $P\in\p3$, the plane selected by 
\,$\omega$
\,is
$a_1(P)x_1+\cdots+a_4(P)x_4=0$,
 (assuming  $P$ is not a singularity of
$\omega$,
\ie not all $
a_i(P)=0$).  The reader is
kindly referred to \cite{CCJ,CLE,Jou}  
for generalities on the subject.

\section{Syzygies}
The condition $ \sum a_ix_i=0$ in (\ref{dist})
ensures the 1-form descends from $\C^4$ to a
section of $H^0({\p3},\Omega^1_{\p3}(m+2))$. It
also tells us that the vector of coefficients,
\ $\uu a:=(a_1,\dots,a_4)$, belongs to the 
  module of syzygies    of $x_1,x_2,x_3,x_4$.
We register for later use the following elementary     
\subsection{\bf Lemma.} 
{\em The syzygies of the regular sequence}
$x_1,\dots ,x_4$ \ {\em are the linear
  combinations} $ \uu a:=(a_1,\dots,a_4)=
\alpha_{12}\kappa_{12}
+\cdots+\alpha_{34}\kappa_{34}
$ 
{\em of the}   \  {\it basic six   Koszul  relations}
$
\kappa_{12}:=(x_2,-x_1,0,0), \dots,
\kappa_{34}:=(0,0,x_4,-x_3) .
$

    \subsection{}\label{syz}
For the case of  1-forms of degree $m=0$, assumed
henceforth, we have
$$\ba c
\omega:=\sum a_idx_i,\ 
 \deg a_i=1,   \ 
\,\sum a_i x_i=0,
\\\na9
\uu a:=
(a_1,\dots,a_4)=\sum_{i<j}{\alpha_{ij}}
{\kappa_{ij}}, \ {\alpha_{ij}\in\C}.
\ea$$
Hence \,$\omega$ \,can be rewritten as
\be\label{omega}
    \omega=\alpha_{12}(x_2dx_1-x_1dx_2)+\cdots+
\alpha_{34}(x_4dx_3-x_3dx_4), \ \alpha_{ij}\in\C.
\ee
It corresponds to an anti-symmetric 4\vez4 matrix,
\bc 
$  \ba c
\omega\leftrightarrow  
\bsm
0&\alpha_{12}&\alpha_{13}&\alpha_{14}\\
-\alpha_{12}&0&\alpha_{23}&\alpha_{24}\\
-\alpha_{13}&-\alpha_{23}&0 &\alpha_{34}\\
-\alpha_{14}&-\alpha_{24}&-\alpha_{34}&0\\
\esm
\leftrightarrow
\text{{bivector}\ in} \wed2  \s1, \ e.g.,
 \\\na9
\kappa_{ij}\lra
x_jdx_i-x_idx_j\lra x_i\wedge x_j .   
\ea
 $
\ec
\subsection{Definition.}\label{contact}
We say
 $\omega$     is a {
\it  1-form of contact}  if
$\det(\alpha_{ij})\nn0$.

In the sequel we identify the 1-form (\ref{omega})
with either the 4\vez4 anti-symmetric matrix
$(\alpha_{ij}) $ or the bivector
\ \,$\sum\alpha_{ij}x_i\wedge x_j, $ as elements of
the projective space $\p5=\p{}(\wed2\s1).$

\subsection{Normal forms}\label{nml}

Linear algebra tells us that the rank of a 4\vez4
nonzero anti-symmetric matrix is either \ 4 \ or
\ 2. Accordingly there are 2 orbits in
\p{}(\wed2\s1) under the induced action of
GL(\s1):
\bi\item
the open orbit,  formed by the anti-symmetric 
4\vez4 matrices
of rank  4, corresponding to contact forms =  
indecomposable bivectors,      and

\item
the closed orbit, \G := 
 Pfaff--Pl\"ucker   quadric,
corresponding to integrable  1-forms
$\leftrightarrow$  
decomposable bivectors. 

\ei

\bc
$\ba c\text{rank 2:}
\bsm\!
\phantom-0&1&0&0\\
\!-1 &0&0&0\\
\! \phantom-0&0&0 &0\\
\!\phantom-0&0&0&0\\
\esm \phantom{rank 2:}\\  \na2
\updownarrow\,
\\\mbox{\normalsize$
\underbrace{
x_2dx_1-x_1dx_2}_{
\text{closed orbit}}$}
\ea
\ \ \  \ \  
\ba c
\text{rank 4:}
\bsm
\phantom-0&1&\phantom-0&0\\
-1 &0&\phantom-0&0\\
\phantom-0&0&\phantom-0&1\\
\phantom-0&0&-1&0\\
\esm\phantom{rank 4:}
\\\na2
\updownarrow\\\mbox{\normalsize$\underbrace{
x_2dx_1-x_1dx_2+x_4dx_3-x_3dx_4}_{
\text{open orbit}}$}
\ea  $
\ec

We recall
these 1-forms define distributions of planes with
nice geometric
interpretations,\,cf.\,\ref{fig1},\,\ref{fig2}. 

$\escala{1.5}\begin{picture}(2,3.35)(.15,  .121)

\ifpdf
\put( .75,1.014) {
\includegraphics[scale=0.71]
{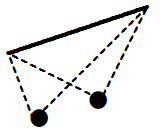}}
\put(4.7, .8) {
\includegraphics[scale=0.39]
{twc.png}}
\fi
\put(.4,.760){rank 2: pencils of planes}
\put(4.5,.760){rank 4: contact distributions}
\end{picture}$

 We'll focus on nested pairs $(\cl D,\varphi)$
consisting of a foliation $\varphi$ of dimension one
and a distribution \cl D of codimension one of
fixed degrees such that  $ \varphi(P)
\subset\cl D(P) \,\forall P\in\p3 $ off singularities.

$\escala{1.5}\begin{picture}(2,2)(.296,2.)
\put(3.4,2.2) {%
\blu{\line(0,1){1.5}%
\line(1,0){1.5}
\line(0,1){1.5}
\put(-1.50,1.5){\line(1,0){1.5}}}
\put(-1.50,1.5){\red{\line(1,-1){1.5}}}
}
\put(2.78,2.1){$\blu{\cl D(P)}$}
\put(4.105,3.204){$\red{\varphi(P)}$}
\put(3.8,2.875){\small$P\ 
\bullet$}
\end{picture}$

\subsection{The space of foliations}
A foliation of dimension one and degree 
$d$  on
\p{{\mb{3}}} is defined by a homogeneous
polynomial vector field
\bc$
\varphi:=p_1\partial_{x_1}+\cdots+p_4\partial_{x_
 {{\mb{4}}}} ,  \,p_i\in\s d
$\ec
up to scalar multiple:
For each point  $P\in\p3$, draw
the line joining 
$P,[p_1(P),\dots,p_4(P)].$ 
The line is undefined at the 
   points where the  2\vez4 matrix
\bc$\bsm x_1&\dots&x_4\\
p_1(x)&\dots&p_4(x)
\esm$
\ec
has rank $<2;$ these are the {\em singular
 points}   of the foliation.  The foliation
remains the same if you add to the
vector field $\varphi$ a multiple of the {radial
 vector field}, ${ \partial_{\rm R}:=
 x_1\partial_{x_1}+ \cdots+x_4\partial_{x_4}} $.
The above matrix changes to
\bc$
\bsm x_1&\dots&x_4\\
p_1(x)+ x_1q&\dots&p_4(x)+x_4q\esm  $;
\ec
getsame Pl\"ucker coordinates (2\vez2 minors).

We often abuse notation and use the same symbol
$\varphi$  for the vector field (assumed\nn0) 
and the foliation of dimenion one it defines.

\subsubsection{\bf Definition}\label{spfolns}
The 
{\em space of foliations}
of dimension one  and  degree  $d$
on \p3 is the  projectivization
\be\label{fd}\left[
\ba l
\p N=\p{}\left(\ds{\frac{
  {\s d\otimes{\s1\us}}}
 {{\s{d-1}\cdot\partial_{\rm R}}}}
\right)=:\F d,  \text{ where}
\\\na{12}{
\s d} = H^0(\p3,\cl O_{\p3}(d)),
\\\na{12}{\s1\us}
=\id{\partial_{x_1},\dots,\partial_{x_4}}=\C^4 ,
\\\na{12}
N+1={4\binom{d+3}3}-{\binom{d+2}3}= (d+4)(d+2)(d+1)/2,
\\\na{15} \left(
\s d\otimes\s1\us\ni\varphi=p_1\partial_{x_1}+\cdots+
p_4\partial_{x_4},\ p_i\in\s d,\text{mod }\partial_{\rm R}
\right).
\ea\right.
\ee
Taking cohomology in the twisted Euler's sequence
$$
\xymatrix
{
\cl O_{\p3}(d-1)\ \ar@{>->}[r]&\cl
O_{\p3}(d)\otimes\s1\us 
\ \ar@{->>}[r]&T\p3(d-1) 
}$$
we have the natural identification
$$\wt\Phi_d:=
H^0(\p3,T\p3(d-1))=\ds{\frac{
  {\s d\otimes{\s1\us}}}
 {{\s{d-1}\cdot\partial_{\rm R}}}}\cdot
$$
Our goal is  to find the dimensions and degrees
of 
certain subvarieties of the space of foliations
$ \F d=\p{}(\wt\Phi_d)$
defined by imposing the
 condition of   tangency to  a varying distribution  
of  planes  in the following sense.

\subsection{Definition}\label{foltgd}
A foliation defined by a vector field 
$\varphi=\sum p_i\partial_{x_i}$ 
is 
{\em tangent}  to  a  
distribution $\omega  =  \sum a_idx_i$ whenever
\bc
$
\omega\cdot\varphi :=  
\sum a_ip_i =  0   $. 
\ec
For fixed $\omega$,
this defines 
a vector \ (resp.\,projective) subspace 
\be\label{perp}\wt
\omega^\perp_d:=\{\varphi\in\wt\Phi_ d\,|\,
\omega\cdot\varphi =  0  \}\   \  
(\rm resp.\ \omega^\perp_d:=\p{}(\wt
\omega^\perp_d)\subset\F d
).
\ee

The game now is to move the 1-form 
$\omega$  
in a suitable family \D\  of distributions, {\it i.e.,} 
pick some closed, irreducible subvariety
\bc$\D\subseteq\p{}(H^0(\p3,\Omega^1(m+2)))$
\ec
  
As customary, we look at the correspondence
$$\ba{cc}
\xymatrix  @R1.32pc@C1pc
{\omega_d^\perp\ar@{->}[d]&\!\!\subset\\
\omega&\!\!\in}
&
\xymatrix  @R1.09pc@C1pc
{
& \wh{\D}_d
:=\{(\omega,\varphi)\in\D\vez\F d\,|\,\omega\cdot\varphi=0\}
\!\!
 \ar@{->}[dl]^{\hskip-.81cm \ba c\vphantom.\\\na{-12.8}
pr_1\ea}
\ar@{->}[dr]_{\hskip .58cm 
\ba c\vphantom.\\\na{-14.8}
pr_2\ea}& 
\\
\D&& \ov{\D}_d\subset\F d.
}\ea$$
The union 
  $$\ov{\D}_d:=pr_2(\wh{\D}_d)=
\ds{\bigcup_{\omega\in\D}}
\omega^\perp_d\ \subset\ \F d
$$
is a first approximation to our object of desire.
Attention must be payed to the {\em variation} of
\ $\dim\omega^\perp_d$ \ as $\omega$ runs in
\D. It may prevent irreducibility of $\wh{\D}_d$.
This is the situation envisaged below, setting
$m=0$ and taking $\bb D=\p5$. Now the fiber
dimensions of ${pr_1}:\wh{\bb D}_d \lar\bb
D$ are given by
\bc
$\dim\omega^\perp_d=\left\{\ba l
2\binom{d+3}3-1 
\text{ for } \omega\in\bb G  \text{ 
(cf.\,(\ref{dimpencil})), whereas it drops to}
\\\na9
 (d+4) (d+2)d /3-1
\text{ for }  \omega\in\p5\setminus\bb G 
\text{ (cf.\,\ref{LdMd})}.
\ea\right.
$
\ec

\noindent
Clearly $\wh {\D}_d$ 
is closed in \D\vez\F d, defined by a bunch of
bilinear equations to be shown in a moment.
  The fiber of $pr_1$ over  $\omega\in\D$  is 
 $\omega^\perp_d$, cf.\,(\ref{perp}).
The image of \ $pr_2$ \ is 

   \centerline{$
\ov{\D}_d=\{ \varphi\, |\,\exists\omega\in\D,
\omega(P)\cdot\varphi(P)=0 ,\forall\,
P\in\p3\}$.
}

\subsubsection{\bf Equations for $\wh{\D}_d$}
We have a natural diagram of maps of vector
bundles over $\D\vez\p3$,
(abusing notation suppressing pullbacks)
$$\ba c
\xymatrix
@R1.65pc@C2.261pc
{\cl O_\D(-1)  \ \ar@{>->}[r]\ar@{->}[dr]&
H^0(\p3,\Omega^1_{\p3}(m+2))\vez\p3
\ \ar@{->>}[d]\\
&\Omega^1_{\p3}(m+2)}
\\\na{-8.94}\ \hskip 1.8621cm
\xymatrix
@R.9961pc{
\ \ar@{= }[d]
\\
\operatorname{Hom}(T\p3,\cl O_{\p3}(m+2))
.}
\ea
$$
The slant arrow 
$\cl
O_\D(-1)\ra \operatorname{Hom}(T\p3,\cl O_{\p3}(m+2))$
yields the {\it universal (twisted) differential form}
\be\label{utd}
T\p3\lar\cl O_{\p3}(m+2)\otimes\cl O_\D(1)
\ee
over \D\vez\p3.
On the other hand, pulling back via $\p3\vez\F d
\ra\F d,$ we have
\bc
$\xymatrix
@R1.5pc@C1.41pc
{\cl O_{\F d}(-1)\ \ar@{>->}[r]\ \ar@{->}[dr]
&H^0(\p3,T{\p3}(d-1))\vez\p3  \ \ar@{->>}[d] &\\&
T{\p3}(d-1)&
\hskip-1.8cm
=\operatorname{Hom}(\cl O_{\p3}(1-d),T\p3 )
}$
\ec
whence the {\it universal vector field} of  degree $d$,
$$\xymatrix@C3.pc
{
\cl O_{\F d}(-1)\otimes\cl O_{\p3}(1-d)
\ \ar@{->}[r]& 
{T\p3} 
\vez\F d 
}
$$
over \p3\vez\F d.
Composing  with the universal differential
form (\ref{utd})
we get a section  
\bc
$\xymatrix
@R1pc@C.91pc
{\cl O_{\p3}(1-d)
\otimes\cl O_{\F d}(-1)
\ \ar@{->}[rr]^{\sigma}&& 
\cl O_\D(1)\otimes\cl O_{\p3}(m+2)
\\&
{T\p3}\ \ar@{<-}[ul]  \ \ar@{->}[ur] &
}$
\ec
over \  $\D\vez\p3\vez\F d$. Twisting by 
$\cl O_\D(-1)\otimes\cl O_{\p3}(d-1)$
we find a diagram
\be\label{sigma1}\ba c
\xymatrix
@R.6522pc
{\hskip-1.8cm
(\sum a_jdx_j)\otimes(\sum p_i\partial_{x_i})
\ \ar@{|->}[rr] && \ \ \sum a_i p_i}
\\
\xymatrix
@R.55  pc
{
\! \ar@{->}[dd]
\cl O_\D(-1)\otimes\cl O_{\F d}(-1)
\ \ar@{->}^{ \hskip  .55cm
\sigma'}[rr]
 &&
\cl O_{\p3}(m+d+1)
\\
& 
\\
\cl O_\D(-1)\otimes{T\p3}(d-1)\ \ar@{->}[uurr]_{\wh\sigma}
&&
}\ea\ee
Taking direct image under 
$\xymatrix
{
pr_{_{\wh3}}:\D\vez\F d\vez\p3\ \ar@{->}[r]& \D\vez\F d
}$, 
we get
$$\ba c\!\!
\xymatrix%
@R.50pc%
@C.152pc
{\cl O_\D(-1)\otimes\cl O_{\F d}(-1)
\ \ar@{->}[dd]
\ \ar@{->}^{\hskip -.69cm\ov\sigma }[rr]&&
\ \ar@{<-}[ddll]
\ \ (pr_3)\ls\left(\cl O_{\p3}(m+d+1)\right)\\
&& 
\hskip 1cm
\sum a_ip_i\in
\stackrel{\mbox{\Large$||$}}
{\s{m+d+1}}
\\  \cl O_\D(-1)\otimes\vphantom{|^{|^|}}
H^0(\p3,T\p3 (d-1))\ni \omega\otimes\varphi
&
&
\!\!\!(\omega=\sum {a_i}dx_i,\, \varphi=\sum{}p_i\partial_{x_i}).
}\\\na{-46.2}
\ifpdf\rotatebox{10} {
\hskip3.0721cm
\xymatrix
@R.379pc
@C.39 pc
{
&
~\hskip1cm
&
\\
 \ar@{|->}[urr]&
&
}}\fi
\ea
$$

\medskip
On the fiber over 
$(\omega=\sum {a_i}dx_i,\  
\varphi=\sum{}p_i\partial_{x_i})
\in\D\vez\F d $
we have \ \ $\ov\sigma(\omega\otimes\varphi)=
\sum {a_i}p_i\in\s{d+m+1}$.  
Set 
$$\left\{\ba l
J^\omega={\id{a_1,\dots,a_4}}, 
\text{ homogeneous ideal of singularities of }\omega,
\\  \na6
(J^\omega)_k=\s k\cap J^\omega,
\text{ subspace of forms of degree $k$ in }J^\omega.
\ea\right.$$ 
For fixed $\omega=\sum {a_i}dx_i
\in\D$, the  vector fields
tangent  to \,$\omega$\,  
form  the  vector subspace as in (\ref{perp}),
\be\label{omegadot}
 \wt
\omega^\perp_d :=\ker\left(\!\!
\ba c
\xymatrix%
@R.4pc@C2.9pc
{H^0(\p3,T\p3(d-1))
\ \ar@{->>}[r]^{\hskip.8cm\omega\cdot
}
&(J^\omega)_{d+m+1}
\\
\sum p_i\partial_{x_i}\ \ \ar@{|->}[r]& \ \ \sum a_ip_i
}
\ea\!\!
\right).
\ee

As said before, we specialize to distributions of
degree $ m=0$.  That's for now 
the sole situation we know
how to control the ranks of $(J^\omega)_d$.
Presently,\bc
$H^0(\p3,\Omega^1(2))=\wed2\s1,\ 
\D\subseteq\p{}(\wed2\s1)=\p5=\{\rm distributions
\,
of\,degree \,0\}
.  
$ 
\ec
Recall 
there are just two orbits (cf.\,\ref{nml}):

\subsubsection{\bf closed orbit}\label{clor}
$^{\mb c}\omega
:=x_2dx_1-x_1dx_2\in\bb G,$ 
pencil of planes.
\\$
\text{We have }
^{\mb c}\omega\cdot
(H^0(\p3,T\p3(d-1)))
= (J^{^{\mb c}\omega})_{d+1}=\id{x_1,x_2}_{d+1}.$
\subsubsection{\bf open orbit}\label{opor}
$^{\mb o}\omega
:=  x_2dx_1  -  x_1dx_2  +  x_4dx_3  -  x_3dx_4,$
\\1-form of contact (\ref{contact}).
Now   $^{\mb o}\omega\cdot
(H^0(\p3,T\p3(d-1)))
=\s{d+1}.$

We look at each at a time, starting with $\bb D=\bb
G.$
\section{Foliations tangent to a pencil of planes}
\noindent
The distributions pictured in \S\ref{fig1}
are given by a 1-form like
$\omega:=
udv-vdu,u,v\in\s1$.
Now  $\omega$ lies in the grassmannian
\G{}=closed orbit in \p5. 
The  rank of the evaluation map \,{$
\omega\mb\cdot
$}  
defined in\,(\ref{omegadot})
is \,$\binom{d+4}3-(d+2)$,
independent of $\omega\in\bb G$:  
\bc
$(J^\omega)_{d+1}
=\id{u,v}\s d=\ker\left(
\s{d+1}\surj{\rm Sym}_{d+1}(\s1/\id{u,v})
\right)
.$
\ec  
Recall the tautological sequence of vector bundles 
of rank 2 over
\G{},
\be\label{RQ}
\xymatrix
{
\cl R\ \ar@{>->}[r]&\s1\vez\G\ \ar@{->>}[r]&\cl Q  
}
\ee
where the fibers $\cl R_{\id{u,v}}=\id{u,v},\ 
\cl Q_{\id{u,v}}=\s1/\id{u,v}.$
It induces the exact sequence
\be\label{clP}
\ba r
\xymatrix%
@C2.963 pc
{\Pi_d\ \ar@{>->}[r]& 
\wed2\cl R\otimes\!
\overbrace{H^0(\p3,T\p3(d-1))}^{
\ba c
(\s d\otimes\s1\us)\big/(\s{d-1}\partial_{\rm R})
\\\na3||\ea
}
\, \ar@{->>}[r]& \cl P_d
}
\\\na9
\xymatrix%
@C1.33 pc
{&( x_2dx_1-x_1dx_2)\otimes\sum 
p_i\partial_{x_i}  \,\ar@{|->}[r]& 
x_2p_1-x_1p_2
}
\ea
\ee
{where }\bc$ 
\cl P_d:=
\ker\left(\G\vez\s{d+1}
\surj
{\rm Sym}_{d+1}\cl Q\right).$
\ec
Tensoring the top row of (\ref{clP}) by $\wed2\cl R\us
=\wed2\cl Q$ we get the exact sequence
\be\label{Pi1}
 \xymatrix
{\Pi'_d:=\Pi_d\otimes\wed2\cl Q\ \ar@{>->}[r]&
 H^0( \p3 ,T\p3(d-1)) \ar@{->>}[r]&\cl P_d 
\otimes\wed2\cl Q.
}
\ee

\subsubsection{\bf Definition}\label{foltgpencil}
The image  $\ov\Pi_d\subset\F d{}$
of the projective subbundle $
\p{}(\Pi'_d)\subset\G{}\vez\F d$
is the {\em variety of foliations
of degree $d$
tangent to some pencil of planes}.
\be\label{dia}
\xymatrix
@C.1pc
{&\ \ar@{->}[dl]_{\mb q }
\p{}(\Pi'_d)\ \ar@{->}[dr]^{\mb p }
&\subset&\G{}\vez\F d  \ar@{->}[dr]^{pr_2} &
\\\bb G&&
\ov\Pi_d & \subset&
\F d.
}
\ee

\subsection{Proposition.}  \label{genbij1}
{\em The  map     \ \ ${\mb p }:
     \p{}(\Pi'_d)\lar\ov\Pi_d\subset\F d$
in the diagram $(\ref{dia})$    
 is generically bijective for   $d  >  1. $
}

\medskip
\noindent
{\em Proof.}
Let $\varphi$ be a general point in
$\ov\Pi_d$. We must show $\varphi$ is tangent to
one and only one distribution $\omega\in\G$.
Say $\varphi=\sum p_i\partial_{x_i},
\omega=x_2dx_1-x_1dx_2$. Now $0=
\omega\cdot\varphi=x_2p_1-x_1p_2$ implies
$p_1=x_1q,p_2=x_2q$ for some $q\in\s{d-1}$, and no
condition is imposed on $p_3,p_4$. If $\varphi$
were tangent to  another $\eta\in\G$, we may
assume  $\eta=x_4dx_3-x_3dx_4$ 
(resp.   $\eta'=x_4dx_2-x_2dx_4$). 
This forces  
$p_3,p_4$ (resp.  $p_2,p_4$)
share a common nonconstant factor
(as $d>1$). 
We have used the fact
that the natural action of GL(\s1) on \G\vez\G\ has
only the orbits of $(\omega,\omega),
(\omega,\eta),
(\omega,\eta')
$ corresponding to the diagonal, the open orbit
and the ``incidence''.
\qed

\section{A polynomial formula}

\subsection{\bf Corollary.}\label{degpen}
{\em For all $d>1$, the degree of the variety}
\ $\ov\Pi_d$ {\em of 
foliations tangent to a 
pencil of planes is given by the top dimensional
Segre class} \ $\mb s_\mb4\Pi'_d$.

\medskip
\noindent
{\em Proof.}
The assertion comes pretty much from Fulton's
construction of Segre \& Chern 
classes,\,\cite[p.\,47]{Ful} taking
(\ref{genbij1}) into account.  
 Setting $h:=$ hyperplane class in $\F d$,
 we have from (\ref{dia})
$$
 \deg\ov\Pi_d=\deg(\ov\Pi_d\cap h^{\dim\ov\Pi_d})
 =\deg \mb p\us h^{\dim\ov\Pi_d}  
 =\deg  \mb q\ls h^{\dim\p{}(\Pi_d)}  
 =\mb s_\mb4\Pi'_d.
$$
\qed

\subsection{Explicit calculation}

The exact sequences (\ref{Pi1}),\,(\ref{clP})
yield 
\bc
$\mb s_\mb4\Pi'_d=\mb c_\mb4(\cl
P_d\otimes\wed2\cl Q
).  $
\ec

Using maple, we may invoke
\cite{kss}{Katz\&Str{\o}mme's Schubert:}

\smallskip

{
\baselineskip 1  pt \parskip 0  pt 

\normalsize

\verb!with(schubert):grass(2,4,q,all);!  

\verb!w2q:=wedge(2,Qq);! \#Pl\"ucker line bdle \verb! O_G(1)!

\verb!Pd:=Symm(d+1,4)-Symm(d+1,Qq);!  

\verb!Pd:=Pd*w2q;!  

\verb!chern(4,Pd);!  

\verb!integral(%);!  

\verb!factor(%);!  

\verb!athus:=unapply(%,d);!  
}



\smallskip
We find the formula for the 
{\em degree of the variety
of foliations of dimension 1 and degree $d$, 
tangent to a (varying)
pencil of planes} in \p3:
\be
\deg\ov\Pi_d=
5\binom{d+4}5 \binom{d+3}3  (d^2+2 d+3)
(d^2+6 d+11)  /108.
\label{athus1}
\ee
Note the degree (=12) in $d$ is thrice the dimension of
the family of distributions at hand. This is expected
by the argument in \S\,\ref{pol} below.

We also register  the calculation of the dimension, 

\be\label{dimpencil}
\left[
\ba r\dim\ov\Pi_d=
\dim\bb G+{\rm rank}\Pi_d-1
=\\
\text{\small (known from (\ref{clP}))}\ \ \ \ \ 
3+
h^0(T\p3(d-1))-{\rm rank}\cl P_d=\\
3+(d+3)(d+2)(d+1)/3=3+2\binom{d+3}{3}
.
\ea\right.
\ee   

\section{Legendrian  vector  fields}

These are polynomial vector fields in \p{3}
(more generaly, \p{{2n+1}})
which are tangent to a 
distribution of contact\,(\ref{opor}).

\subsection{Example.}
The   vector field 
$
 x_1^2\partial_{x_1}+x_2x_1\partial_{x_2}
+x_3x_4\partial_{x_3}+x_4^2\partial_{x_4}   $ 
is  tangent to
$^{\mb o}\omega
:=  x_2dx_1  -  x_1dx_2  +  x_4dx_3  -  x_3dx_4$.

\subsection{Universal anti-symmetric  map}
We'll need
the universal anti-symmetric  map  $\alpha$ 
as defined at the bottom of the diagram below:
\be\label{uas}
 \ba c\xymatrix%
@R1.49pc@C1pc
{
\left(   \s1\us \otimes\wed2\s1\right.
\ \ar@{->}[rr]&& 
\left.\vphantom{\wed2\s1}\ \ \s{1}\right)_{\p 5}
\\
\s1\us \otimes\cl O_{\p5}(-1)\vphantom{\ba c|\\|\ea}
\ \ar@{>->}[u]
\ \ar@{->}[rr]^{\hskip.8cm \alpha}
& &\,\s{1}\otimes\cl O_{\p5}\!\!\!\hskip-.51cm
\ar@{= }[u]
\vphantom{\ba c|\\-\ea}
}

\\\na{-1}
\xymatrix
@C.81pc
{
\partial_{x_i}\otimes \sum u_j\wedge{}v_j\ 
\ar@{|->}[rr]& &\ 
\sum(\partial_{x_i} u_j)\cdot v_j
-(\partial_{x_i} v_j)\cdot u_j
}
\ea\ee

The rank of $\alpha$   drops from  4 to 2 over 
  \bb G = Pfaff--Pl\"ucker quadric in \p5=
\p{}(\wed2\s1). Recall \bb G\ is the locus of
decomposable 1-forms (bivectors)
$ {\omega=udv-vdu(\leftrightarrow{}u\wedge{}v)}$. 

Tensoring (\ref{uas}) with 
$\s d=H^0(\p3,\cl O_{\p3}(d))$, we get the diagram
\be\label{uasd}
\xymatrix
@R1pc{
\s d\otimes\s1\us \otimes{\cl O_{\p5}(-1)}
\ \ar@{->}[rr] ^{\hskip.51cm1_{\s d}\otimes\alpha}
\ar@{->}[rrdd]_{\tau_d}& & \s d\otimes\s{1}{}_{
  |\p5}\ \ar@{->}[dd]
\\& &
\\
\s{d-1}\cdot\partial_{\rm R} \otimes{\cl
O_{\p5}(-1)}\vphantom{\ba c|\\|\ea}
\ar@{>->}[uu] 
\ \ar@{->}[rr]^{\hskip1cm 0}& 
&\s{d+1}{}_{ |\p5}
}
\ee
  The map     ${\tau_d}$    passes  to 
  the  quotient    and  induces
\be\label{ovtaud}
\ba c
 \left(
\ds{\frac{\s{d}\otimes \s1\us}
{\s{d-1}\cdot  \partial_{\rm
R}}}
\right) \otimes{\cl O_{\p5}(-1)}
\xymatrix
@C3.1pc
{
\ \ar@{->}[r]^{\hskip-.65 cm\ov\tau_d}&  \ \,  
\s{d+1}{}_{|\p5}  
}
\\\na6
\ \ \ \ 
(\varphi\mod\partial_{\rm R})\otimes{\omega}
\ 
\xymatrix
@C1.9pc
         {\ \ar@{|->}[r]& 
          }
\ \
{a_1}p_1+\cdots+{a_4}p_4
\ea  
\ee
with notation as in display\,(\ref{fd}).

\subsubsection{\bf Remarks.}\label{dirim}
We register the map $\ov\tau_d$ in display
 (\ref{ovtaud}) arises via direct image of 
$$
T\p3(d-1)\otimes\cl O_{\p5}(-1)\stackrel{\wh\tau_d}
\lar\cl O_{\p3}(d+1)
\ \text{ by }\ 
pr_1:\p5\vez\p3\ra\p5
$$
as in (\ref{sigma1}). A crucial observation is the fact
that the image sheaf of $\ov\tau_d$ is the direct image
of the image sheaf of $\wh\tau_d,\,d>>0$:   kill
the appropriate $H^1$. The map
${\ov\tau_d}$ 
 is  surjective  off  the  Pfaffian 
 $=\bb G\subset
\p{}(\wed2\s1)$ because so is ${\tau_d}$. 

 Let  
\be\left\{
\ba l  \cl L_d:=
 \ker\left({\ov\tau_d}\otimes\cl O_{\p5}(1)\right)  
\,\subset\,    H^0(\p3,T\p3(d-1));
\\\na9
\cl M_d:={\rm image}\left({\ov\tau_d}\otimes
\cl O_{\p5}(1)\right)  
\,\subset\,   \s{d+1}\otimes\cl O_{\p5}(1).

\ea\right.
\ee

$  \cl L_d$ and $\cl M_d$ are
 vector  bundles over the 
open subset  \ $\p5 \setminus \bb G.$

\subsection{Lemma.} \label{LdMd}  
$\cl L_d,\cl M_d$
{\em both
  extend  as vector  bundles over} \p5 
{\em and fit into
the exact sequence}
\\\centerline{$
\xymatrix
@C1.85pc
{
\cl L_d \ \ar@{>->}[r]& 
H^0({\p3},T\p3(d\!-\!1))_{| \p5}    \ 
\ar@{->>}[rr]^{\!\!
\ov\tau_d\otimes1_{\cl O_{\p5}(1)}}
&&
\cl M_d\subset
\s{d+1}\otimes\cl O_{\p5}(1).
}$}

\smallskip\noindent
{\em Moreover, 
\\\mb{(i)} the inclusion $
\cl M_d\subset
\s{d+1}\otimes\cl O_{\p5}(1)$ is an equality over 
 $\p5 \setminus \bb G$ and 
\\\mb{(ii)} we have {\rm rank} $\cl L_d
=(d+4) (d+2) d /3
$.
}

\medskip\noindent
{\it Proof.} This is because
the rank of  $\ov\tau_d$  drops precisely
 along a   hypersurface. Indeed, looking at the diagram
(\ref{uasd}), notice the rightmost vertical map is
 surjective, whereas the top horizontal one drops rank
precisely along the Pfaffian.
Now use the following (inspired by 
Raynaud-Gruson's \cite[(5.4.3)]{ray})
\subsubsection{\bf Claim} {\em
Let $\theta:\cl E_1\ra\cl E_2$ be a map of locally free
sheaves over a variety $X$. Let $r$ denote the
generic rank of $\theta$. 
Let \cl J be the image sheaf of $\wed r\cl E_1
\otimes\wed r\cl E_2\ve\ra\cl O_X$ induced by $\theta$.
Let $\G^r_{\cl E_1}$ be the
Grassmann bundle of rank $r$ quotients of ${\cl E_1}$.
There is a rational section $
X\dashrightarrow\G^r_{\cl E_1}
$ 
induced by $\theta$. 
It extends to a morphism if \cl J is
invertible. In this case, the image sheaf of $\theta$
is locally free.
}

\subsubsection{\bf Proof of  Claim}
Let $X'\ra X$ be the blowup of \cl J.
Arguing with the Pl\"ucker embedding
$\G^r_{\cl E_1}\subset\p{}(\wed{r}\cl E_1)$, we
see that the composition $X'\ra X\dashrightarrow
\G^r_{\cl E_1}
$ is a morphism. Since $X'\ra X$ is an isomorphism for
\cl J invertible, we are done.
\qed

\medskip
Over the open orbit, 
 $\cl L_d$ \  is  the  set  of  pairs
  $(\omega,\varphi)$
such that  the vector field  $\varphi$
is \ Legendrian w.r.t. the form of 
contact $\omega$.

The image  $\cl M_d$ is a locally free subsheaf 
of $\s{d+1}\otimes\cl O_{\p5}(1)$ (of same rank), 
though not locally
split.     

\subsubsection{\bf Definition.}
Notation as in (\ref{fd}),
the image $\bb L_d  \subset  \F d$   
of  the  projective  subbundle
$$
\xymatrix
@C1pc
{&
\p{}(\cl L_d)\ \ar@{->}[dr]^{\mb p }
&\subset&\p5  \vez\F d  \ar@{->}[dr]^{pr_2} &
\\\p5\ \ar@{<-}[ur]_{pr_1}
&&
\bb L  _d & \subset&
\F d 
}
$$
is \ the variety of {\it Legendrian   \ foliations}.
 
\subsection{Proposition.}  \label{genbij}
{\em The  map     \ \ ${\mb p }:
     \p{}( \mathcal{L}_d)\lar
     \mathbb{L}_d\subset\F d$
 is generically bijective for   $d  >  1. $
}

\medskip
\noindent
{\em Proof.}
Let $\varphi$ be a general point in
$\mathbb{L}_d$. 
We must show $\varphi$ is tangent to
one and only one distribution $\omega\in\p5\setminus\G$.
We have from \cite[Theorem 7.2]{CJM} that
 $\varphi$ corresponds to a
decomposable  2-form 
$\alpha \wedge \omega$, where $\alpha$ is a
polynomial   1-form of degree $d$.   
If $\varphi$
were tangent to  another $\eta\in\p5\setminus\G$,  
then  
$\alpha \wedge \omega \wedge \eta =0$.  By
\cite[Theorem 5.1]{CJM}  the singular points of  $
\omega \wedge \eta$ consist of two skew lines, in
particular has  codimension $2$. So  Saito's   
division Lemma \cite{Saito} applies and
we get  that  $\alpha
=  f\omega + g \eta$, 
 where  $f,g\in\s{d-1}$. Therefore,
$\alpha \wedge \omega=g \cdot (\eta \wedge
 \omega)$. We arrive at a contradiction, since
 $\deg(g)=d-1>0$ and $\alpha \wedge \omega$ has no
 zeros of codimension one.  
\qed

\subsubsection{\bf Remark}
Alan  Muniz has communicated to the authors that
the decomposability  result
\cite[Theorem\,7.2]{CJM}  also follows from Saito's    
division Lemma \cite{Saito}.   In fact, let
$\theta$ be a homogeneous
polynomial 2-form on $\mathbb{C}^4$
  which induces a Legendrian
foliation on \p3. 
Consider the  linear contact form $\omega$ such
that $\theta \wedge \omega=0$. Since
$\{\omega=0\}=\{0\}\subset  \mathbb{C}^4$,  then 
 by Saito's    division Lemma   there is a
 polynomial  1-form $\alpha$  such that
 $\theta=\alpha \wedge \omega$.   

\smallskip
\subsubsection{\bf Corollary.} 
{\em The degree of} $\bb L_d\ (d>1)$
{\em  is given by the  
Segre class}  $\mb s_\mb5\cl L_d=\mb c_\mb5\cl M_d
$.

\medskip{\em Proof.} Same as for tangency to
pencils of planes, cf.\,\ref{degpen}.\qed

\subsection{Calculation}
We proceed to the actual calculation. Recall the
exact  sequence      (\ref{LdMd}).
Since the  middle term  is trivial,
it implies  
$$
\rm Segre\cl L_d=(Chern\cl L_d)\inv=Chern\cl M_d    .
$$
The degree of the latter class can be found
{\em using Bott's formula in the  equivariant
context} as we learn from \cite{Meurer}.

As kindly enticed by the referee, a practical
summary of the main ingredients  
is included in the sequel for 
the reader's convenience. We are given a smooth,
projective variety $X$ endowed with an action of
the torus $\bb T=\G_m$. Let $\cl E\stackrel\pi\lar X$
be an equivariant vector bundle. So
there is an action $
\G_m\vez\cl E\stackrel\psi\lar\cl E$ yielding 
a natural commutative diagram
$$\ba{ccc}
\G_m\vez\cl E&\stackrel\psi\lar&\cl E\\
\downarrow&&\ \,\downarrow\pi\\\G_m\vez
X&\lar& X.
\ea
$$
If $P\in X$ is a point fixed by the action, it
follows that $\G_m$ acts on the fiber $\cl E_P$.
Any such action splits $\cl E_P=\oplus \cl E_P^{\chi}$,
a direct sum of eigenspaces with character
$\chi$. This means $\G_m\vez\cl E_P\ni(t,v)\mapsto
t\cdot v=\sum \chi(t)v_\chi$. Each $\chi(t)=t^{w_\chi}$
for some integers $w_\chi$, called weights.
Bott's formula expresses integrals of polynomials
in the Chern classes of equivariant vector bundles
in terms of the corresponding classes in
equivariant cohomology. With the simplifying
assumption that the set of points in $X$  left fixed
by the action is finite, the equivariant cohomology
classes are just symmetric functions on the 
 weights. Next we make this explicit for $c_5\cl M_2$.

  This requires the  description of the fibers  of
  \ $\cl M_d$
over the \mb6 fixed points of \p{}(\wed2\s1)
under an
action of  \bb C\us\  induced from a  natural action 
on \s1: 
$x_i\mapsto t^{w_i}x_i$\,, \,where the weights
$w_i\in\bb N$
will be chosen appropriately.
  The fixed points are  the ``canonical''
  anti-symmetric matrices     of rank  2,
namely

\smallskip
{\small$
\bsm
0 & 1 & 0& 0\\
\!\!-1\,\,& 0 & 0& 0\\
0 & 0 & 0& 0\\
0 & 0 & 0& 0\\
\esm  ,
\bsm
0 & 0 & 1& 0\\
0 & 0 & 0& 0\\
\!\!-1\,\,& 0 & 0& 0\\
0 & 0 & 0& 0\\
\esm\!,
\bsm
0 & 0 & 0& 1\\
0& 0 & 0& 0\\
0 & 0 & 0& 0\\
\!\!-1\,\,& 0 & 0& 0\\
\esm\!,
\bsm
0 & 0 & 0& 0\\
0 & 0 & 1& 0\\
0 & \!\!-1\,\,& 0& 0\\
0 & 0& 0& 0\\
\esm,
\bsm
0 & 0 & 0& 0\\
0 & 0 & 0& 1\\
0 & 0 & 0& 0\\
0& \!\!-1\,\,& 0& 0\\
\esm\!,
\bsm
0 & 0 & 0& 0\\
0 & 0 & 0& 0\\
0 & 0 & 0& 1\\
0 & 0 &\!\!-1\,\,& 0\\
\esm\!.  $}%

\smallskip

Let us work out the fiber of $\cl M_d$ say over 
the last matrix in the above list.  That matrix
 corresponds  to the 1-form
\be\label{w0}
\omega_0:=x_4dx_3-x_3dx_4. \ee
We perturb $\omega_0$  to full rank, moving 
away from the Pfaffian. 
Define   
\be\label{omegat}
\omega_{t}:=x_4dx_3-x_3dx_4+ t (x_1dx_2-x_2dx_1).
\ee
The fiber of $\cl M_d$ over $\omega_0$
is found by evaluating $\omega_t$
on vector fields of degree  $d$   and taking
$\lim_{t\ra0}$.
In the sequel, set for simplicity, $d=\mb2.$
The calculations reported below were performed
using {\sc singular}\,\cite{sing}; script
available in \cite{arXiv}.

We  construct a basis for
$H^0(\p3,T\p3(\mb2-1))$
consisting of eigenvectors for the
given
\C\us-action.  
Consider the subspace
\\\centerline{ $
{\s2^\circ}:
=\{\varphi:=\Sigma p_i\partial_{x_i}
\in\s2\otimes\s1\us
\st p_i\in\s{2},\Sigma\partial_{x_i}p_i=0
\}$}
consisting of
vector fields  with  zero  divergence.   

We have
\ $\left(\s1\cdot
\partial_{\rm R}\right)\oplus{\s2^\circ}=\s2\otimes\s1\us,  $
\ whence
$H^0(\p3,T\p3(1))
\raise-2pt\hbox{$
\ba c{\wt~}\\\na{-16.}-\ea$}
{\s2^\circ}$.

We want a basis for the space
$$
\s2^\circ=\ker\left(\!\!\ba c
\xymatrix
@C3.1pc
{
\s2\otimes\s1\us\ \ar@{->>}[r]^{\ \,
\rm tr(jac)}& 
\s1}\\
\xymatrix
{
\Sigma p_i\otimes\partial_{x_i}\ \ar@{|->}[r]& 
\Sigma
\partial_{x_i}p_i
}
\ea\!\!
\right).  
$$

It's formed, when $d=2$, by  
\mb{ 40-4 } vectors like (cf.\,\cite[L245]{arXiv})
\be\label{base}
\ba{cccr}
[x_2^2,0,0,0],&\dots&[-x_1x_3,0,0,x_3x_4],&\dots\\
|| &  \dots&||&\dots
\\
x_2^2\partial_{x_1},&\dots&-x_1x_3\partial_{x_1}+
x_3x_4\partial_{x_4},&\dots
\\
{2w_2-w_1},&\dots&{w_1\!+\!w_3\!-\!w_1}
={w_3\!+\!w_4\!-\!w_4},&\dots
\ea 
\ee
{The bottom row} in the above display
lists the weights of the 
elements in the basis.
   Compute next the matrix of the map
$\s2^\circ\lar\s3$ defined via contraction by the 
perturbed 1-form $\omega_t$
(\ref{omegat}):

\centerline{$\ba l\\
\omega_t[x_2^2,0,0,0]=-tx_2^3,  
\dots,\\\dots
\\
\omega_t[-x_1x_3,0,0,x_3x_4]=-x_3^2x_4- tx_2x_1x_3 
,\dots.
\ea  $}

\noindent
The sought for fiber of the image $\cl M:=\cl M_2$
is computed via saturation:  

\bi\item
perform elementary row operations on the matrix of
$\omega_t ; $

\item divide each row by its gcd;
 
\item last but not least,  set   $t=0.$
 
 \ei
We get a matrix with the correct rank \mb{20}.   
The columns of pivots form a basis for the fiber of
\cl M$_2$\ at the fixed point $\omega_0$ cf.(\ref{w0}).

Say for instance the basic vector
 $[x_2^2,0,0,0]
=x_2^2{\partial_{x_1}}
 $ 
corresponds to a column with a pivot.
Then the fiber of \cl M$_2$\ 
acquires an eigenspace which is 
a direct summand with weight
 $2w_2{-w_1}.
$ That's all we care to collect and eventually
obtain numbers feeding into Bott's formula.

For $\omega_t$ as in (\ref{omegat}),
letting $t\ra0$, 
we find the fiber of \cl M$_2$\ 
is the direct sum of
20 (=rank) eigenspaces with corresponding weights
$$\ba c
2w_1-w_2,w_1,w_2,-w_1+2w_2,-w_3+2w_4,w_4,w_3,2w_3-w_4,
\\
w_2-w_3+w_4,w_2,w_2+w_3-w_4,2w_2-w_3,2w_2-w_4,w_1-w_3+w_4,\\w_1,w_1+w_3-w_4,w_1+w_2-w_3,w_1+w_2-w_4,2w_1-w_3,2w_1-w_4
\ea$$
Plugging in numerical values for the $w_i$, {\it e.g.,}
0,2,7,10,    we find the list of 
20 weights
\\\centerline{$
-2,0,2,4,13,10,7,4,5,2,-1,-3,-6,3,0,-3,-5,-8,-7,-10.
$}

The equivariant (top=5=DIM) Chern class at the
chosen fixed point $\omega_0$
is 
the value of the 5th elementary symmetric
function of those weights,
{\it i.e.,} 
the coefficient of $t^5$ in the 
product
$(1-2t)\cdot{}(1+0t)\cdot{}(1+2t)\cdots(1-7t)
(1-10t)
$,
 to wit,\bc$
\fbox{\bf105534}$.\ec
Repeating the calculation for each of the 6 fixed
points, 
Bott's   formula reads
\bc
$
\deg{\rm Chern}\cl M_2
=\sum_P c^{\rm equiv}_5\cl M_P\big/
c_5^{\rm equiv}T_P\p5  =$\ec

$ \frac{833800359}{42000}
\!-\!\frac{38740434}{1500}\!-\!
\frac{  4199874 }{336}\!+\!
\frac{  7716777 }{  336 }\!-\!\frac{  3398841 }{1500}
\!-\!\frac{\fbox{\bf105534}}{  42000 }
 =\mb{2224}.$


\medskip
 This is the value of the new
{\bf Athus \pol\ }\,(\ref{athusbis})
for  $d=2.$

\section{\bf Another  
\pol\ formula}\label{pol}

The calculation for higher degrees are better left 
for a script, cf. appendix of arXiv version \cite{arXiv}
of this article.
We find by interpolation,
\be\label{athusbis}\ba c
 {{\binom{d+2}4}}\cdot  (d^3+9d^2+14d+24)\cdot\\\na8
\left(\vphantom{|^|_|}d^8+34d^7+475d^6+3430d^5+13480d^4+\right.
\\\na8
\left.29872d^3
   +45444d^2+44856d+29808\vphantom{|^|_|}\right) /  (2 ^53^55).
\ea\ee

%


%
%

As
we shall argue, the
formula for the degree is a polynomial in $d$ of degree
$\leq15$, so it  suffices to compute for the
16 values $d\in\{2,3,\dots,17\}$. Polynomiality is
deduced from the fact that $\cl M_d$ is a direct
image in \p5 
of a sheaf on \p5\vez\p3 (cf.\,\ref{dirim})
whose Chern character
is a \pol\ in $d$ of degree $\leq3$. Using
Grothendieck-Riemann-Roch \cite[(5)\,p.\,283]{Ful}
it follows that the
Chern character of $\cl M_d$ is a \pol\ in $d$ of
degree $\leq3$. Hence the top Chern class  $c_5\cl M_d$
is a \pol\ with degree $\leq15$ as asserted.

\section{Final remarks}

The case of  Legendrian vector fields on \p5
can be
treated similarly, except for the need to perform
an initial blowup of  \p{14}, the space of
  anti-symmetric 6\vez6 matrices, 
 along the grassmannian of lines 
in \p5.  Ditto for higher odd dimensional.
However, the computational load as of now
seems unfeasible: It took several many days to 
work out the very first few values.
\be\label
{P5}
\!\!
\left\{\ba{ll}
2:&310143368560\\  
3:&241876654493880936\\ 
4:&11354802732747615971781\\
5:&96189569307604075178197866\\
6:&249521135387730096977922116592\\
7:&269205473509858802653259153925591
\\
8:&146758500496340866823126747040755640
\ea\right.\ee
Code available as an appendix to
the
arXiv version of this note\,\cite{arXiv}.

\end{document}